# MASS EXTINCTIONS: AN ALTERNATIVE TO THE ALLEE EFFECT

By Rinaldo B. Schinazi

*University of Colorado*

We introduce a spatial stochastic process on the lattice $\mathbf{Z}^d$ to model mass extinctions. Each site of the lattice may host a flock of up to $N$ individuals. Each individual may give birth to a new individual at the same site at rate $\phi$ until the maximum of $N$ individuals has been reached at the site. Once the flock reaches $N$ individuals, then, and only then, it starts giving birth on each of the $2d$ neighboring sites at rate $\lambda(N)$. Finally, disaster strikes at rate 1, that is, the whole flock disappears. Our model shows that, at least in theory, there is a critical maximum flock size above which a species is certain to disappear and below which it may survive.

**1. Introduction and results.** It seems that the main mass extinction theory proposed in ecology, for species that reproduce sexually, is the so-called Allee effect: if the density of a certain species is driven sufficiently low, then encountering mates of the opposite sex becomes unlikely and the population is driven to extinction even if left alone by predators or disease; see Stephens and Sutherland (1999). While the Allee effect seems suitable to explain extinction of animals living by themselves or in small flocks, it does not look suitable to explain the extinction or near extinction of animals such as passenger pigeons which apparently remained in large flocks almost to their end; see Austin (1983).

We propose a mathematical model that, at least in theory, shows that animals living in large flocks are more susceptible to mass extinctions than animals living in small flocks. More precisely, we will show that if the maximum flock size is above a certain threshold, then the population is certain to become extinct, while if the maximum flock size is below the threshold, there is a strictly positive probability that the population will survive.

Our model is a spatial stochastic model on the lattice $\mathbf{Z}^d$, typically $d = 2$. Each site of the lattice may host a flock of up to $N$ individuals. Each individual may give birth to a new individual at the same site at rate $\phi$ until the maximum of $N$ individuals has been reached at the site. Once the flock reaches $N$ individuals, then, and only then, it starts giving birth on each of the $2d$ neighboring sites at rate $\lambda(N)$. This rule is supposed to mimic the fact that individuals like to stay in a flock and will give birth outside the flock only when the flock attains the maximum number $N$ that a site may support. Finally, disaster strikes at rate 1, that is, the whole flock disappears. This rule mimics an encounter with greedy hunters or a new disease. Both disasters seem to have stricken the American buffalo and the passenger pigeon.

We now write the above description mathematically. Each site $x$ of $\mathbf{Z}^d$ may be in one of the states: $0, 1, 2, \ldots, N$ and this state is the size of the flock at $x$. The model is a continuous-time Markov process that we denote by $\eta_t$. Let $n_N(x, \eta_t)$ be the number of neighbors of site $x$, among its $2d$ nearest neighbors, that are in state $N$ at time $t$.

Assume that the model is in configuration $\eta$; then the state at a given site $x$ evolves as follows:

$$i \to i+1 \qquad \text{at rate } i\phi + \lambda(N) n_N(x, \eta) \text{ for } 0 \leq i \leq N-1,$$
$$i \to 0 \qquad \text{at rate 1 for } 1 \leq i \leq N.$$

Schinazi (2002) has introduced a model related to the present one for a different question.

We will be interested in extinctions in two different senses. We say that finite populations die out if, starting from any finite population, there is a finite random time after which all sites are empty. We say that infinite populations die out if, starting from any infinite population, for any given site there is a finite random time after which the site will be empty forever.

The model in the special case $N = 1$ is well known and is called the contact process [see Liggett (1999)]. For the contact process, we know that there exists a critical value $\lambda_c$ (that depends on the dimension $d$ of the lattice) such that the population dies out (in the two senses defined above) if and only if $\lambda \leq \lambda_c$.

We now state our main result.

THEOREM 1. *Consider the model with parameters $N$, $\lambda(N)$ and $\phi$ and let*

$$m \equiv 2d\lambda(N) \prod_{i=1}^{N-1} \left(1 - \frac{1}{1 + i\phi + 2d\lambda(N)}\right) \qquad \text{for } N \geq 2,$$
$$m \equiv 2d\lambda(1) \qquad \qquad \qquad \qquad \qquad \text{for } N = 1.$$

(a) *Assume $m \leq 1$; then any finite population dies out.*
(b) *Assume $m < 1$; then any infinite population dies out.*

Next we show that if $\lambda(N)$ does not grow too rapidly with $N$, then the population dies out for large $N$.

COROLLARY 1. *Assume that $\lambda(N)$ and $\phi > 0$ are such that*

$$\lim_{N \to \infty} \frac{\lambda(N)}{N^{1/(1+\phi)}} = 0;$$

*then any finite or infinite population dies out when $N$ is large enough.*

A little calculus is enough to prove Corollary 1 so we skip the proof. Note also that the conclusion of Corollary 1 also holds when $\phi = 0$ provided there exists $a < 1$ such that

$$\lim_{N \to \infty} \frac{\lambda(N)}{N^a} = 0.$$

Our main application is the following easy consequence of Corollary 1.

COROLLARY 2. *Assume $\lambda(N) \equiv \lambda$ and that $\lambda > \lambda_c$ (the critical value of the contact process), $\phi > 0$. Then, there is a critical positive integer $N_c(\lambda, \phi)$ such that any finite population dies out for $N > N_c$ and survives for $N < N_c$. The same is true for infinite populations for a possibly different critical $N_c$.*

PROOF. According to Corollary 1, finite and infinite populations die out for $N$ large enough for constant $\lambda$. On the other hand, if $\lambda > \lambda_c$ we know that the model with $N = 1$ has a positive probability of not becoming extinct. Finally, as a consequence of the construction provided in Section 2 we will see that for constant $\lambda$, the smaller the $N$ the more likely it is for the population to survive. Putting together these three facts, we get the existence of the critical value $N_c$. This completes the proof of Corollary 2. □

We believe that the critical value $N_c$ is the same for finite and infinite populations but this is still unproved.

Next we show that a low internal birth rate $\phi$ may be compensated by a large external birth rate $\lambda(N)$ but that the converse is not true.

THEOREM 2. (a) *For all $\phi \geq 0$, $N \geq 1$ and all initial configurations starting with at least one site in state $N$, there is a positive probability, for finite or infinite populations, not to die out, provided $\lambda(N)$ is large enough.*

(b) *If $\lambda(N) < \lambda_c$, then finite and infinite populations die out for all $\phi \in [0; \infty]$.*

**2. Construction of the process and proof of Theorem 1.** We now give an explicit graphical construction for the process $\eta_t$. Let $\|\cdot\|$ denote the Euclidean norm. Consider a collection of independent Poisson processes: $\{L^{x,y}, F^{x,i}, D^x : x, y \in \mathbf{Z}^d, \|x-y\| = 1, 1 \leq i \leq N-1\}$. For $x$ and $y$ in $\mathbf{Z}^d$ such that $\|x - y\| = 1$, let the intensity of $L^{x,y}$ be $\lambda(N)$. For $x$ in $\mathbf{Z}^d$ and an integer $i$ between 1 and $N-1$, let the intensity of $F^{x,i}$ be $i\phi$. Finally, for $x$ in $\mathbf{Z}^d$, let 1 be the intensity of $D^x$. The graphical construction takes place in the space-time region $\mathbf{Z}^d \times (0, \infty)$. At an arrival time of $L^{x,y}$ ($\|x - y\| = 1$), if site $x$ is in state $N$ and there are $N-1$ or fewer individuals at site $y$, then we add an individual at $y$. At an arrival time of $F^{x,i}$ and if $x$ is in state $i$, $1 \leq i \leq N-1$, then we change the state of $x$ to $i+1$. Finally, at an arrival time of $D^x$ we put $x$ in state 0. In this way we obtain a version of our spatial stochastic process with the precribed rates. For more on graphical constructions, see, for instance, Durrett (1995).

Assume $N_1 < N_2$ and $\lambda(N) \equiv \lambda$. Using the graphical construction above, construct the model, $\eta_{2,t}$, for parameters $(\lambda, N_2, \phi)$. We can also construct the model $\eta_{1,t}$ with parameters $(\lambda, N_1, \phi)$ in the same probability space by using the same Poisson processes $L^{x,y}$, $D^x$ and by using the processes $F^{x,i}$ only for $i \leq N_1 - 1$. Start $\eta_{1,t}$ and $\eta_{2,t}$ with a single individual at the origin. Both processes are in the same configuration until they reach state $N_1$ at the origin. At this point in time, the flock at the origin for $\eta_{1,t}$ starts giving birth to individuals in neighboring sites while the flock at the origin for $\eta_{2,t}$ continues increasing internally. Since the death rates and the external birth rate $\lambda$ are the same for both processes, it is easy to check the following. No transition can break the inequalities

$$\min(\eta_{2,t}(x), N_1) \leq \eta_{1,t}(x).$$

In this sense, the lower the $N$ the more spread out the population and the more likely it is to survive. Note that this coupling works only for constant $\lambda$. The fact that the model is more likely to survive for small $N$ is all that was missing to the proof of Corollary 2.

One can see that some attempted births will not occur because the site on which the attempted birth takes place has already reached the maximum size $N$. This creates dependence between the size of the offspring of different individuals. Because of this lack of independence, explicit probability computations seem impossible. In order to prove Theorem 1, we introduce a branching-like process for which explicit computations are possible and that dominates, in a certain sense, our process $\eta_t$.

We now describe the new process informally. It may be constructed in the same way as $\eta_t$ by using appropriate Poisson processes. For a formal construction of a similar process, see Pemantle and Stacey (2001). While for $\eta_t$ there is a maximum of one flock per site (with a maximum size of $N$), for the new branching-like process, that we denote by $b_t$, there is no limit on the number of flocks per site but each flock is again limited to $N$ individuals maximum. For $b_t$, as for $\eta_t$ and with the same rate $\lambda(N)n_N$, each new flock is started by a birth from one of its neighbors. However, for $b_t$, once a flock is started it grows only through internal births. That is, a flock that has started does not receive births from neighbors. We take the internal birth rate for a flock in $b_t$ to go from $i$ to $i+1$ to be $i\phi + 2d\lambda(N)$. Note that this rate is the maximum growth rate for a flock in $\eta_t$ (that rate is achieved only if all its $2d$ neighbors are in state $N$). Each flock of $b_t$ that reaches size $N$ starts giving birth to individuals in neighboring sites at rate $\lambda(N)$. Each of these births starts a new flock, since there is no bound on the number of flocks per site for $b_t$. Finally, each flock dies, independently of everything else, at rate 1.

We now give a more mathematical description of the process $b_t$. Each site $x$ of $\mathbf{Z}^d$ is in state 0 (empty) or in state $(r, i_1, \ldots, i_r)$, where $r \geq 1$ represents the number of flocks at site $x$ and $i_1, \ldots, i_r$ represent the number of individuals (between 1 and $N$) of each flock. Let $n_N(x, b_t)$ be the number of flocks of size $N$ in the neighborhood of site $x$. The transition rates for $b_t$ at a site $x$ are given by

$$(r, i_1, \ldots, i_r) \to (r+1, i_1, \ldots, i_r, 1) \quad \text{at rate } \lambda(N)n_N(x, b_t),$$

$$(r, i_1, \ldots, i_r) \to (r, i_1, \ldots, i_{j-1}, i_j + 1, i_{j+1}, \ldots, i_r) \quad \text{at rate } i_j\phi + 2d\lambda(N)$$
$$\text{for } 1 \leq j \leq r \text{ and if } i_j \leq N-1,$$

$$(r, i_1, \ldots, i_r) \to (r-1, i_1, \ldots, i_{j-1}, i_{j+1}, \ldots, i_r) \quad \text{at rate } 1$$
$$\text{for } 1 \leq j \leq r \text{ and } r \geq 2,$$

$$0 \to (1, 1) \quad \text{at rate } \lambda(N)n_N(x, b_t),$$

$$(1, i_1) \to 0 \quad \text{at rate } 1.$$

Note that birth rates are higher for $b_t$ than for $\eta_t$, that death rates are the same and that all attempted births actually occur for $b_t$ while they may or may not occur for $\eta_t$. Techniques such as in Liggett [(1985), Theorem 1.5 in Chapter III] can be used to construct the processes $b_t$ and $\eta_t$ in the same probability space in such a way that, if they start with the same initial configuration, if there is a flock of size $i$ on a site $x$ for $\eta_t$, then there is at least one flock of size at least $i$ for $b_t$ on the same site $x$.

We start by proving Theorem 1(a).

Consider the process $b_t$ starting with a single individual at the origin of $\mathbf{Z}^d$. We call such an individual, who is the first individual of a flock, a founder. We are going to compute the expected number of founders a given founder gives birth to. Let $A$ be the event: "the founder's flock will eventually reach the maximum size $N$." In order to reach $N$, the flock must add one individual at a time, $N-1$ times before getting wiped out. Using properties of the exponential distribution one gets

$$P(A) = \prod_{i=1}^{N-1} \frac{i\phi + 2d\lambda(N)}{1 + i\phi + 2d\lambda(N)}.$$

Let $X$ be the number of founders given birth to by a single founder. In order to give birth to $k$ founders, the founder must first start a flock that will reach $N$ and then the flock must give birth $k$ times before disappearing. Again by using properties of the exponential distribution we get for $k \geq 1$

$$P(X = k) = \left(\frac{2d\lambda(N)}{2d\lambda(N) + 1}\right)^k \frac{1}{2d\lambda(N) + 1} P(A).$$

Therefore, the expected number of founders is

$$E(X) = 2d\lambda(N) P(A).$$

Note that $E(X) = m$, where $m$ has been defined in the statement of Theorem 1.

Let the first founder be the zeroth generation and let $Z_0 = 1$. This first founder gives birth to a random number of founders, before dying, and these form the first generation. Denote their number by $Z_1$. More generally, let $n \geq 1$; if $Z_{n-1} = 0$, then $Z_n = 0$; if $Z_{n-1} \geq 1$, then $Z_n$ is the total number of founders the $Z_{n-1}$ founders, of the $(n-1)$st generation, give birth to before dying. It is clear that the process $Z_n$ is a Galton–Watson process and it dies out if and only if

$$m \leq 1.$$

Note that if $Z_n$ becomes extinct, so does $b_t$ and therefore $\eta_t$. It is easy to see that the same is true if we start with any finite number of individuals instead of 1. This completes the proof of Theorem 1(a).

We now prove Theorem 1(b).

Let $\bar{\eta}_t$ be the process $\eta_t$ starting with $N$ individuals per site. Using the graphical construction one can construct $\bar{\eta}_t$ and $\eta_t$ so that $\bar{\eta}_t(x) \geq \eta_t(x)$ for every site $x$, all times $t > 0$ and any initial configuration $\eta_0$. Therefore, it is enough to show that $\bar{\eta}_t$ dies out. That is, for any site $x$ in $\mathbf{Z}^d$ there is a time $T$ such that $\bar{\eta}_t(x) = 0$ for all $t > T$. We will actually prove this claim for the process $b_t$, starting with one flock of size $N$ at each site.

Note that if there is at least one individual at site $x$ at time $t$ for the process $b_t$, then it must be the case that this individual is the descendent of an individual who was on some site $y$ at time 0. Let $Z_n(y)$ be the number of founders [as defined in the proof of Theorem 1(a)] of the $n$th generation of the process started at $y$ with one flock of $N$ individuals. For an individual at $y$ to be the ancestor of an individual at $x$, the process $Z_n(y)$ must have survived at least $\|x - y\|$ generations. This is so because each generation gives birth on nearest-neighbor sites only. So $n$th-generation founders are at distance $n$ or less from $y$. Since there is one flock of $N$ individuals at time 0 at $y$, the expected number of founders in generation 1 is $2d\lambda(N)$. From generation 1 onwards $Z_n(y)$ is a Galton–Watson process with mean offspring $m$. Thus, let $n = [\|x-y\|] + 1$, where $[a]$ is the integer part of $a$, and using that $m < 1$ we get

$$P(Z_n(y) \geq 1) \leq E(Z_n(y)) = 2d\lambda(N)m^{n-1} \leq 2d\lambda(N)m^{\|x-y\|-1}.$$

Therefore,

$$\sum_{y \in \mathbf{Z}^d} P(Z_n(y) \geq 1) \leq \sum_{y \in \mathbf{Z}^d} \frac{2d\lambda(N)}{m} m^{\|x-y\|} < \infty.$$

The Borel–Cantelli lemma implies that almost surely there is an integer $\ell$ such that if $\|y - x\| > \ell$, then $y$ cannot be an ancestor of $x$. On the other hand, according to Theorem 1(a), any finite population dies out. Thus, the population which was initially on sites $y$ such that $\|y - x\| \leq \ell$ is dead after a finite random time $T$. This shows that site $x$ remains empty after $T$.

This completes the proof of Theorem 1(b).

**3. Proof of Theorem 2.** Theorem 2(a) can be proved in a pretty standard way so we will only sketch its proof. We deal with the case $d = 1$: it is easy to see that if the process survives in $d = 1$, then it will survive in any other dimension. Let

$$B = (-4L, 4L) \times [0, T].$$

Assume that each site of $[-L, L]$ is in state $N$. Consider the process restricted to the space-time region $B$; that is, assume that there are no births from outside $B$ into $B$.

Let $\varepsilon > 0$; it is easy to see that we can pick $T$ (depending on $\varepsilon$ and $L$) so that the probability that there are no deaths in the space-time box $B$ is at least $1 - \varepsilon/2$.

Note that even if $\phi = 0$, the flock at $L$ can fill up the site $L + 1$. This can happen provided the state at $L$ is $N$. Once the site $L+1$ is in state $N$ it can start filling the site in $L+2$ and so on. More formally, we define the process $r_t$ (as rightmost site) as $r_0 = (L, N)$, for which the first coordinate indicates the site position and the second coordinate indicates the state of that site. The evolution rules for $r_t$ are given by

$(i, N) \to (i+1, 1) \quad$ at rate $\lambda(N)$ for $L \leq i \leq 4L - 1$,

$(i, j) \to (i, j+1) \quad$ at rate $\lambda(N)$ for $L+1 \leq i \leq 4L-1$ and $1 \leq j \leq N-1$.

In the absence of deaths in $B$, all the sites between $-L$ and $3L$ will be in state $N$ by time $T$ provided $r_t$ reaches $(3L, N)$ by time $T$. Ignoring the internal births, that is, assuming that $\phi = 0$, only delays the filling process. It is easy to see that counting all the successive transitions for $r_t$ gives a Poisson process $R_t$ with rate $\lambda(N)$. We have that $r_T = (3L, N)$ if and only if $R_T$ is at least $(2L+1)N$. This has probability at least $1 - \varepsilon/4$ provided $\lambda(N)$ (that depends on $N$, $L$, $T$ and $\varepsilon > 0$) is large enough. The same may be done to show that, in the absence of deaths, with probability at least $1 - \varepsilon/4$ all the sites between $-3L$ and $L$ will be filled by time $T$. Therefore, one sees that with probability at least $1 - \varepsilon$ one block of sites, in $[-L, L]$, in state $N$ gives birth to two blocks of sites, in $[-3L, -L]$ and in $[L, 3L]$, in the same state. Moreover, this is true for the process restricted to $B$ and uniformly on all possible states of the boundary of $B$.

Now, well-known techniques [see, e.g., Theorem 4.4 in Durrett (1995)] allow us to compare the process $\eta_t$ to a (very) supercritical oriented percolation on

$$\mathcal{L} = \{(m, n) \in \mathbf{Z}^2 : m + n \text{ is even}\}.$$

This comparison implies survival for finite and infinite populations.

Note that if at least one site is in state $N$, then there is a positive probability, even if $\phi = 0$, to get a block of $2L+1$ sites in state $N$ and we may start the construction above. This completes the sketch of the proof of Theorem 2(a).

We now turn to the proof of Theorem 2(b). It is essentially the same proof as the one of Theorem 1(b) in Schinazi (2002). Since it is short we include it.

Consider the model with $\phi = \infty$. In this case, as soon as there is one individual at a site it immediately fills to $N$. So each site has only two possible states: 0 and $N$. The transition rates are given by

$$0 \to N \quad \text{at rate } \lambda(N) n_N(x, \eta),$$

$$N \to 0 \quad \text{at rate } 1.$$

Therefore, the model above is a contact process with birth rate $\lambda(N)$. If $\lambda(N) < \lambda_c$, this contact process dies out.

Using the graphical construction, it is easy to see that the model with $\phi = \infty$ can be coupled to a model with any finite $\phi$ in such a way that, starting from the same configuration, the model with $\phi = \infty$ has more individuals, site per site, than the model with finite $\phi$. Since the model with $\phi = \infty$ dies out, so does the model with a finite $\phi$. This completes the proof of Theorem 2(b).

**Acknowledgment.** We thank Enrique Andjel for several helpful suggestions.

DEPARTMENT OF MATHEMATICS
UNIVERSITY OF COLORADO
COLORADO SPRINGS, COLORADO 80933-7150
USA
E-MAIL: schinazi@math.uccs.edu
URL: http://math.uccs.edu/~schinazi/